\documentclass{amsart}
\usepackage{longtable}
\usepackage{amssymb,amsthm}
\usepackage{hyperref}

\newtheorem{theorem}{Theorem}[]

\newtheorem{corollary}[theorem]{Corollary}

\theoremstyle{definition}
\newtheorem{example}[theorem]{Example}
\newtheorem{question}[theorem]{Question}

\begin{document}
\title[On the proportion of derangements and on suborbits]{On the proportion of derangements and on suborbits in finite transitive groups}

\author[M.~Barbieri]{Marco Barbieri}
\address{Dipartimento di Matematica ``Felice Casorati",\newline%
University of Pavia, Via Ferrata 5, 27100 Pavia, Italy} 
\email{marco.barbieri07@universitadipavia.it}
	
\author[P.~Spiga]{Pablo Spiga}
\address{Pablo Spiga, Dipartimento di Matematica Pura e Applicata,\newline%
University of Milano-Bicocca, Via Cozzi 55, 20126 Milano, Italy} 
\email{pablo.spiga@unimib.it}

\thanks{The authors are members of the GNSAGA INdAM and PRIN ``Group theory and its applications'' research group and kindly acknowledge their support.}
\dedicatory{With deep appreciation for Professor G\'abor Korchmáros's insights on derangements, we dedicate this paper with admiration.}

\maketitle

\begin{abstract}
	We find a lower bound on the proportion of derangements in a finite transitive group that depends on the minimal nontrivial subdegree. As a consequence, we prove that, if $\Gamma$ is a $G$-vertex-transitive digraph of valency $d\ge 1$, then the proportion of derangements in $G$ is greater than $1/2d$. 
\end{abstract}

\section{Introduction}
\label{sec:introduction}

Let $G$ be a finite transitive permutation group with domain $\Omega$. A \emph{derangement of $G$} is a permutation $g$ without any fixed points, that is, for any $\alpha \in \Omega$ we have that $\alpha$ and $\alpha^g$ are distinct points of $\Omega$.

C.~Jordan has first noticed in \cite{Jordan} that every nontrivial finite transitive permutation group contains a derangement. This observation has far reaching consequences ranging from number theory to topology (see, for instance, \cite{FeinKantorSchacher,Serre}). For further information, we refer to the extensive review in \cite[Chapter~1]{BurnessGiudici}.

An outstanding property of derangements in transitive groups is their abundance. We define the \emph{proportion of derangements of $G$} as the ratio
\[ \delta (G) = \frac{|F_0(G)|}{|G|} , \]
where $F_0(G)$ denotes the set of derangements in $G$. We can give two impressive examples of lower bounds for $\delta(G)$. In \cite{CameronKovacsNewmanPraeger}, P.~J.~Cameron, L.~G.~Kov\'acs, M.~F.~Newman and C.~E.~Praeger have proven that, if $G$ is a $p$-group, then
\[ \delta(G) > \frac{p-1}{p+1} . \]
On the other end of the spectrum, there exists an absolute constant $\varepsilon > 0$ such that, for every simple group $T$ endowed with a transitive action, $\delta(T) \ge \varepsilon$. This result has been obtained by J.~E.~Fulman and R.~M.~Guralnick by an extensive study of primitive actions of finite simple groups \cite{FulmanGuralnick1,FulmanGuralnick2,FulmanGuralnick3,FulmanGuralnick4}.

In this paper, we prove a lower bound on $\delta(G)$ depending on the minimal nontrivial subdegree of the transitive group $G$ (see Section~\ref{sec:mainResult}), and, as a consequence, we obtain a lower bound for the proportion of derangements in the automorphism group of a symmetric graph depending on the valency of the graph alone (see Section~\ref{sec:graphCorollary}).

\section{A new bound for $\delta(G)$}
\label{sec:mainResult}

Before continuing, we need to introduce some notion. Fix a point $\alpha\in \Omega$. We let
\[ G_\alpha = \{g\in G \mid \alpha^g=\alpha\} \]
denote the \emph{stabilizer of $\alpha$ in $G$}.
The orbits of $G_\alpha$ acting on $\Omega$ are called the \emph{suborbits of $G$ with respect to $\alpha$}. Let
\[ O_1, O_2, \ldots, O_r \]
be the list of all suborbits of $G$. The integer $r$ is called the \emph{permutational rank of $G$}, while the cardinalities 
\[ |O_1|=d_1, \ldots, |O_r|=d_r \]
%of the orbits of $G_\alpha$
are called the \emph{subdegrees of $G$}. Without loss of generality, we can suppose that $O_1=\{\alpha\}$: this suborbit is the \emph{trivial suborbit}, and its cardinality $d_1$ is the \emph{trivial subdegree}.\footnote{Observe that $d_1,\ldots,d_r$ are not necessarily distinct, hence $G$ might have more than one subdegree equal to $1$.}
We define
$$d_G=\min\{d_2,\ldots,d_r\}$$
and we call $d_G$ the \emph{minimal nontrivial subdegree of $G$}.\footnote{We remark that $d_G$ is well-defined except when $r=1$, that is, $\Omega=\{\alpha\}$. In this degenerate case, we let $d_G=\infty$.}

In \cite{CameronCohen}, P.~J.~Cameron and A.~M.~Cohen have shown that, for every transitive groups $G$ with permutational rank $r$,
\[ \delta(G) \ge \frac{r-1}{|\Omega|}. \]
Furthermore, this bound is sharp: indeed, the equality is achieved if and only if $G$ is a Frobenius group.\footnote{In~\cite{MR3383274}, R.~M.~Guralnick, I.~M.~Isaacs and the second author have proven an upper bound on $\delta(G)$ depending on $r$ only.}

Our main result is heavily inspired by their proof, with a shift in focus from the permutational rank to the minimal nontrivial subdegree.

\begin{theorem}\label{thrm:main}
	Let $G$ be a finite transitive permutation group whose minimal nontrivial subdegree is $d_G$. Then
	\[ \delta(G) \ge \frac{1}{2d_G} + \frac{n-2}{2|G|} . \]
	Equality is attained if and only if $G$ is a Frobenius group.
\end{theorem}

Theorem~\ref{thrm:main} is proved in Section~\ref{sec:proof}, and our bound is compared with the Cameron--Cohen bound in Section~\ref{sec:comparison}. We are not able to determine a criterion for choosing, \emph{a priori}, which bound is optimal.

\section{A graph-theoretic corollary}
\label{sec:graphCorollary}

What prompted us to study whether $\delta(G)$ could be bounded from below by a function of the minimal nontrivial subdegree $d_G$ is a graph theoretic problem. In the \emph{10th PhD Summer School in Discrete Mathematics} in Rogla, G.~Korchm\'aros has asked the first author whether, for every graph $\Gamma$ of valency $d$ and transitive automorphism group $\mathrm{Aut}(\Gamma)$, the proportion of derangements in $\mathrm{Aut}(\Gamma)$ is bounded away from zero by a function of $d$ alone. We would like to thank him for inspiring us with this fascinating problem.

The connection between our group theoretic perspective and the graph theoretic problem is achieved through the notion of \emph{orbital digraphs}. Every permutation group $G$ can be endowed with a family of digraphs whose vertex set is the domain $\Omega$ and whose arcs are $G$-orbits of a given pair $(\alpha,\beta) \in \Omega \times \Omega$. We recall that the orbital digraphs and the suborbits of $G$ are in one-to-one correspondence, and the image of the suborbit $O_i$ under this map is a regular digraph of valency $d_i$. Moreover, every digraph admitting a group of automorphism $G$ is a union of orbital digraphs of $G$.\footnote{We say that $G$ is a \emph{group of automorphism} for a digraph $\Gamma$ if $G$ is a subgroup of $\mathrm{Aut}(\Gamma)$.} We refer to \cite[Section~3.2]{DixonMortimer} for a more exhaustive account.

Theorem~\ref{thrm:main} gives at once an affirmative answer to G.~Korchm\'aros's question.
\begin{corollary}
	Let $\Gamma$ be a finite digraph, and let $G$ be a group of automorphisms of $\Gamma$. If $G$ is transitive, and $\Gamma$ has valency $d$, then 
	\[ \delta(G) \ge \frac{1}{2d} . \]
\end{corollary}

\section{Proof of Theorem~\ref{thrm:main}}
\label{sec:proof}

We start by establishing some notation for the proof. Let $G$ be a nonidentity finite transitive group of degree $n$ with domain $\Omega$ and let $\alpha\in \Omega$. Given $g\in G$, we let
\[ \mathrm{Fix}(g)=\{\omega\in \Omega\mid\omega^g=\omega\} , \]
and, given $i\in \{0,\ldots,n\}$, we let
\[ F_i(G)=\{g\in G\mid |\mathrm{Fix}(g)|=i\} . \]
In particular, $F_0(G)$ is the set of all derangements of $G$. Our aim to show that
\[ \delta(G) = \frac{|F_0(G)|}{|G|} \ge \frac{1}{2d_G} + \frac{n-2}{2|G|} , \]
and that the equality is attained if and only if $G$ is a Frobenius group. 

Since the sets $F_i(G)$ partition $G$, we get 
\begin{align}\label{eq:1}
	|G|&=\sum_{i= 0}^n|F_i(G)| .
\end{align}
Moreover, from the Orbit Counting Lemma (see \cite[Theorem~1.7A]{DixonMortimer}), we have
\begin{align}\label{eq:2}
	|G|&=\sum_{i=0}^ni|F_i(G)| .
\end{align}
Observe that $F_n(G)$ consists of the identity of $G$. By subtracting Equation~\eqref{eq:1} from Equation~\eqref{eq:2}, we deduce
\begin{align*}
	|F_0(G)|&=\sum_{i= 1}^n(i-1)|F_i(G)|
	\\&=\sum_{i=1}^{n-1}(i-1)|F_i(G)|+n-1
	\\&\ge\sum_{i= 2}^{n-1}|F_i(G)|+n-1
	\\&=\sum_{i=2}^n|F_i(G)|+n-2
	\\&=|G|-|F_0(G)|-|F_1(G)|+n-2.
\end{align*} 
Therefore
\begin{align}\label{eq:3}
	|F_0(G)|&\ge\frac{|G|}{2}-\frac{|F_1(G)|}{2}+\frac{n-2}{2}.
\end{align}

Observe that that the sets $F_1(G_\omega)$, as $ \omega$ runs through the elements of $\Omega$, are pairwise disjoint and cover the whole of $F_1(G)$. This means that 
\[\{F_1(G_\omega)\mid\omega\in \Omega\}\]
is a partition of $F_1(G)$ and hence, for every fixed $\alpha \in \Omega$,
\begin{align}\label{eq:7}
	|F_1(G)|&=\sum_{\omega\in \Omega}|F_1(G_\omega)|=|\Omega||F_1(G_\alpha)|=\frac{|G|}{|G_\alpha|}|F_1(G_\alpha)|.
\end{align}

Now, let $\beta\in\Omega\setminus\{\alpha\}$ such that $|\beta^{G_\alpha}|=d_G$, and let $G_{\alpha,\beta} = G_\alpha \cap G_\beta$. Since $G_{\alpha,\beta}\cap F_1(G_\alpha)=\emptyset$, we get $F_1(G_\alpha)\subseteq G_\alpha\setminus G_{\alpha,\beta}.$ Therefore,
\begin{align}\label{eq:8}|F_1(G_\alpha)|\le |G_\alpha|-|G_{\alpha,\beta}|=|G_\alpha|\left(1-\frac{|G_{\alpha,\beta}|}{|G_\alpha|}\right)=|G_\alpha|\left(1-\frac{1}{d_G}\right).\end{align}

Finally, combining Equations~\eqref{eq:3},~\eqref{eq:7} and~\eqref{eq:8}, we get
$$	\delta(G)=\frac{|F_0(G)|}{|G|}\ge\frac{1}{2}-\frac{1}{2}\left(1-\frac{1}{d_G}\right)+\frac{n-2}{|G|}=\frac{1}{2d_G}+\frac{n-2}{2|G|}.$$

Moreover, if the equality is attained, then, from Equation~\eqref{eq:3}, we deduce that
\[ G = F_0(G) \cup F_1(G) \cup F_n(G) , \]
that is, $G$ is a Frobenius group. Conversely, it is readily seen that any Frobenius group attains our bound.

\section{Comparison}
\label{sec:comparison}

In this section, we compare the bound obtained in Theorem~\ref{thrm:main} with the Cameron--Cohen bound from \cite{CameronCohen}.

For every $\alpha \in \Omega$, we note that
\begin{align*}
	(|G_\alpha|+d_G) n &= (|G_\alpha|+d_G)\sum\limits_{i=1}^r d_i 
	\\&= |G_\alpha| + d_G + |G_\alpha| \left( \sum\limits_{i=2}^r (d_i-d_G)+d_G(r-1) \right)
	\\&\qquad+ d_G\left( \sum\limits_{i=2}^r (d_i - |G_\alpha|) + |G_\alpha|(r-1) \right)
	\\&= |G_\alpha| + d_G + 2d_G|G_\alpha|(r-1)
	\\&\qquad+ d_G|G_\alpha|\sum\limits_{i=2}^r \left( \dfrac{d_i}{d_G} + \dfrac{d_i}{|G_\alpha|} - 2 \right) .
\end{align*}
By substituting this equality in the difference of the two bounds, we obtain
\begin{align*}
	\frac{1}{2d_G} + \frac{n-2}{2|G_\alpha|n} - \frac{r-1}{n}
	&= \frac{|G_\alpha|n + d_G(n-2)}{2d_G|G_\alpha|n} - \frac{r-1}{n}
	\\&= \frac{(|G_\alpha|+d_G)n - 2d_G}{2d_G|G_\alpha|n} - \frac{r-1}{n}
	\\&= \frac{1}{2n} \left(\frac{1}{d_G} - \frac{1}{|G_\alpha|} + \sum\limits_{i=2}^r \left( \dfrac{d_i}{d_G} + \dfrac{d_i}{|G_\alpha|} - 2 \right) \right) .
\end{align*}
Therefore, the sign of
\begin{align} \label{eq:c2}
	\frac{1}{d_G} - \frac{1}{|G_\alpha|} + \sum\limits_{i=2}^r \left( \dfrac{d_i}{d_G} + \dfrac{d_i}{|G_\alpha|} - 2 \right)
\end{align}
determines which bound gives the best estimate for $\delta(G)$.
In particular, the Cameron--Cohen bound is better when Equation~\eqref{eq:c2} is negative, while our bound is stronger otherwise. We remark that the sign of Equation~\eqref{eq:c2} depends on the distribution of the nontrivial subdegrees of $G$: a predominance of subdegrees proximate to $|G_\alpha|$ results in a positive expression, whereas a prevalence of subdegrees closer to the minimal nontrivial subdegree leads to a negative sign. 

We conclude this section by giving two examples of infinite families of transitive permutation groups: two in which the bound in Theorem~\ref{thrm:main} is stronger, and one in which the Cameron--Cohen bound is better.

\begin{example}
	Let $\mathrm{PSL}_2(p)$ be the projective special linear group of dimension $2$ over the field with $p$ elements, $p$ prime, and suppose that $p \equiv 43 \pmod{120}$. This condition guarantees that the alternating group $\mathrm{Alt}(4)$ is a maximal subgroup of $\mathrm{PSL}_2(p)$ and simplies some computations. We consider the primitive action of $\mathrm{PSL}_2(p)$ with stabilizer $\mathrm{Alt}(4)$. The possible subdegrees of this action are $1, 4, 6, 12$. Let us denote by $\mu_i$ the number of suborbits having cardinality $i$. A computation shows that 
	\begin{align*}
		n &= \frac{p(p^2-1)}{24} ,\\
		\mu_1 &= 1 ,\\
		\mu_4 &= \frac{p-4}{3} ,\\
		\mu_6 &= \frac{p-3}{8} ,\\
		\mu_{12} &= \frac{1}{12}(n - 1 - 4\mu_4 - 6\mu_6) \simeq \frac{p^3}{24 \cdot 12} .
	\end{align*}	
	Moreover,
	$ r= 1+\mu_4+\mu_6+\mu_{12}$ and $n = 1 + 4\mu_4+ 6\mu_6+ 12\mu_{12} .$
	Therefore,
	\[ \lim_{p\to\infty} \frac{r-1}{n} = \frac{1}{12} < \frac{1}{8}=\frac{1}{2d_{\mathrm{PSL}_2(p)}} .\]
	Since most subdegrees are equal to the cardinality of a point stabilizer, as expected, our bound is stronger, for sufficiently large primes $p$.
\end{example}

\begin{example}
	In~\cite{Spiga2024}, the second author of this paper gives remarkable examples of transitive permutation groups where most points lie in a suborbit of cardinality $2$: examples of this type are relevant for the enumeration of vertex-transitive graphs of given valency. In these examples, $|G_\alpha|=4$ and hence the subdegrees of $G$ are $1$, $2$ or $4$. Let $\mu_i$ be the number of subdegrees of $G$ having cardinality $i$. From~\cite{Spiga2024}, we deduce that $\mu_1=n/6$, $\mu_2=n/3$ and $\mu_4=n/24$ and hence
	\[r=\mu_1+\mu_2+\mu_4=\frac{13n}{24}.\]
	Therefore,
	\[\frac{1}{2d_{G}}+\frac{n-2}{2|G|}=\frac{1}{2}+\frac{n-2}{8n}\quad\hbox{and}\quad
	\frac{r-1}{n}=\frac{13n -24}{24n}.\]
	In particular, the bound in Theorem~\ref{thrm:main} is stronger than the Cameron--Cohen bound.
\end{example}

\begin{example}
	Let $G$ be a non-Frobenius $2$-transitive permutation group of degree $n$. The Cameron--Cohen bound is $1/n$, while the bound in Theorem~\ref{thrm:main} is
	\[ \frac{1}{2d_G} + \frac{n-2}{2|G|} < \frac{1}{2(n-1)} + \frac{n-2}{2n (n-1)} = \frac{1}{n} .\]	
\end{example}

We conclude this note with a question.
\begin{question}
	Is there a function $f: \mathbb{N} \to \mathbb{N}$ such that, for every permutation group $G$ of degree $n$, minimal nontrivial subdegree $d$ and rank $r$, if $n \ge f(d)$, then
	\[ \frac{r-1}{n} \le \frac{1}{2d}+\frac{n-2}{2|G|} \,?\]
\end{question}
In essence, our inquiry revolves around determining whether our bound exhibits asymptotic superiority over the Cameron--Cohen bounds when $d$ remains fixed.

\bibliographystyle{plain}
\bibliography{bibDerangemntsBoundedValency}
\end{document}